\documentclass[graybox]{svmult}

\usepackage{mathptmx}       
\usepackage{helvet}         
\usepackage{courier}        
\usepackage{makeidx}         
\usepackage{graphicx}        
\usepackage{multicol}        
\usepackage[bottom]{footmisc}

\usepackage{cite}
\usepackage{fancyvrb}        

\usepackage{amssymb}
\usepackage{amsmath}

\begin{document}

\title*{Latency Exploitation in Wavelet-based Multirate Circuit Simulation}
\author{Kai Bittner 
\and
Hans Georg Brachtendorf 
}
\institute{Kai Bittner, Hans Georg Brachtendorf \at University of Applied Sciences Upper Austria, 
         Softwarepark 11, 4232 Hagenberg, Austria, 
\email{Kai.Bittner@fh-hagenberg.at, Hans-Georg.Brachtendorf@fh-hagenberg.at}}

%
%
\maketitle

\abstract*{
The simulation of radio frequency (RF) circuits is one of the severest problems in 
Design Automation: the information signal or envelope is modulated by a carrier signal
with a center frequency typically in the GHz range. Due to Nyquist's sampling theorem the
time steps in conventional transient analysis are prohibitively small. A technique
to overcome Nyquist's bottleneck is the multirate method which reformulates the
ordinary circuit's differential algebraic equations (DAEs) as a system of partial DAEs (PDAEs).
In this paper further improvements of the wavelet multirate circuit simulation technique are presented.
In the new algorithm we use different grids for the approximation of the solution on
different circuit parts, exploiting latency. 
In particular, for circuits comprising latent parts the grids can be much sparser,
which results in the reduction of the overall problem size and leads to a faster simulation.}

\abstract{The simulation of radio frequency (RF) circuits is one of the severest problems in 
Design Automation: the information signal or envelope is modulated by a carrier signal
with a center frequency typically in the GHz range. Due to Nyquist's sampling theorem the
time steps in conventional transient analysis are prohibitively small. A technique
to overcome Nyquist's bottleneck is the multirate method which reformulates the
ordinary circuit's differential algebraic equations (DAEs) as a system of partial DAEs (PDAEs).
In this paper further improvements of the wavelet multirate circuit simulation technique are presented.
In the new algorithm we use different grids for the approximation of the solution on
different circuit parts, exploiting latency. 
In particular, for circuits with latencies the grids can be much sparser,
which results in the reduction of the overall problem size and leads to a faster simulation.}

\section{Introduction}

In simulation of RF circuits one faces waveforms with a spectrum centered
around a center frequency, which is typically in the GHz range for
modern communication standards. Due to the Nyquist's theorem the waveforms
must be discretized with a sampling rate, which is at least twice as high
as the highest relevant frequency in the spectrum. Classical transient solvers
which solve the initial value problem (IVP) show unacceptably long run times.
To overcome this bottleneck envelope methods based on a reformulation
of the ordinary DAEs by partial DAEs, known als multirate PDAE have been
developed \cite{Bra94,BWL+96,NL96,Roy97,PG02,Bra97,Bra2001,BiBra12b}. However,
despite this tremendous progress, the run time is often prohibitively
long for circuits such as PLLs. In this paper improvements based on
latency exploitation are proposed, which utilize some specific properties
of (sub-)circuits of RF circuitry and properties of the multirate PDAE. 

For speeding up conventional transient analysis, several attempts have
been made for exploiting latency and other specific properties of
circuit DAEs. An excellent overview of these methods can be found,
e.g.,\ in \cite{Ogro94}.
In \cite{New79} a relaxation method denoted as timing analysis has
been presented based essentially on one Gauss-Seidel (GS) iteration per time step.
This method has an emphasis on CMOS circuits where the diagonal part
of the Jacobian matrix is dominant. Alternatively, the 
waveform relaxation (WR), e.g.,\ \cite{LRSV82,FMS95}, has attracted  attention
for several decades. Here the circuit is divided into sub-circuits wherein
the coupling of these sub-circuits is relatively week. Each sub-circuit is
simulated for a time period while the remaining sub-circuits are idle or latent.
The method is repeated until convergence is achieved.
WR may be interpreted as a block Gauss-Seidel for a time period. It has been
developed for CMOS circuits, too. In \cite{SA01} the latency insertion
method (LIM) has been proposed, which has its origin in electromagnetic
field simulation. Essentially, the discretization grid for the unknown currents
and voltages are shifted. The technique is advantageous when delays
of interconnects are dominant. The time steps however are limited similar
to the CFL condition.
Node tearing, often with latency exploitation, has been
reported, e.g.,\  in \cite{CBH+91,SD85,RSH79,RH76}. In \cite{SD85,RH76} 
the sub-circuits
are allowed to have separate integration step sizes reflecting their
activity level.

For all the cited methods a careful partitioning and/or time step control
is required to achieve convergence. The method for the latency exploitation
considered in this paper, which is
based both on the multirate PDAE and spline-wavelet technique, has none of these
restrictions.

\section{The multirate circuit simulation problem \label{MULTI_RATE}}

We consider circuit equations in the charge/flux oriented
modified nodal analysis (MNA) formulation, which yields a mathematical model
in the form of a system of differential-algebraic equations
(DAEs):
\begin{equation}
  \label{eq_MNA_charge}
  \tfrac{d}{dt}q\big(x(t)\big)
        + g\big(x(t)\big) = s(t).
\end{equation}
For RF circuits the circuit DAE (\ref{eq_MNA_charge}) exhibits multirate behavior, i.e.,
(most) of the signal waveforms have a bandpass spectrum, where the spectrum is centered
around a center frequency, which is typically in the GHz range for state of the art
mobile phone standards. The time steps employing conventional solvers for ordinary
DAEs must be kept sufficiently small to avoid aliasing of the numerical solution.
The run time is therefore prohibitively long. One method to
overcome this bottleneck  reformulates the underlying ordinary
DAEs by a system of partial DAEs \cite{Bra94,BWL+96}. Several modifications of this
method have been proposed \cite{NL96,Roy97,PG02,Bra97,Bra2001,BiBra12a}. 

To separate different time scales the problem is
reformulated as a multirate PDAE, i.e.,
\begin{equation}\label{multirate}
\Big(\tfrac{\partial}{\partial \tau} 
+\omega(\tau)\,\tfrac{\partial}{\partial t} \Big) q\big(\hat{x}(\tau,t)\big)
+g\big(\hat{x}(\tau,t)\big)=\hat{s}\big(\tau,t\big)\ 
\end{equation}
with mixed initial-boundary conditions $x(0,t)=X_0(t)$ and $x(\tau,t)=x(\tau,t+P)$.
A solution of the original circuit equations can be found along certain characteristic lines \cite{BiBra12b}. 

\begin{figure}
\centering
\includegraphics[width=\columnwidth]{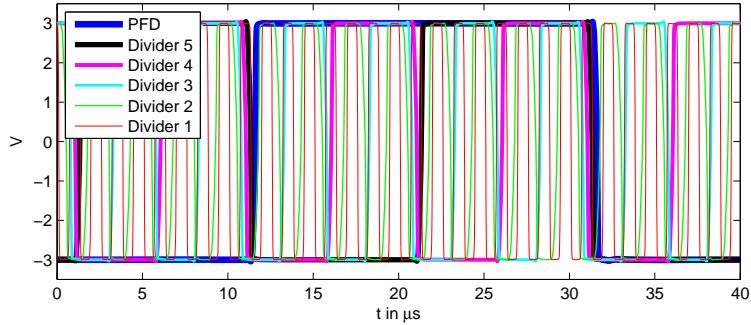}
\caption{Several signals in a frequency divider chain as part of a PLL.}
\label{bittner:fig_1}       
\end{figure}

Discretization with respect to $\tau$ (Rothe's method) using a linear multistep
method results in a periodic boundary value problem in $t$
of the form
\begin{align}\label{envelope}
\omega_k\,\tfrac{d}{dt}q_k\big(X_k(t)) + &f_k(X_k,t) = 0,\\
\nonumber X_k(t) = &X_k(t+P),
\end{align}
where $X_k(t)$ is the approximation of $\hat{x}(\tau_k, t)$ for the $k$-th time step $\tau_k$ (cf.\ \cite{BiBra12b,BiBra14b}).
The periodic boundary value problem (\ref{envelope}) can be solved by several methods, as 
Shooting, Finite Differences, Harmonic Balance, etc. Here, we consider the spline wavelet based method
introduced by the authors in \cite{BiBra14b}, following ideas from \cite{BiDau10a,BiDau10b}. 
One problem of traditional methods is that all signals in the circuit are discretized
over the same grid. This can pose a problem if different signal shapes are present in the circuit, which may be approximated 
more efficiently if individual grids are used for each of the signals. As an example we consider a chain of 5 frequency dividers
(as part of a PLL). In each step the frequency is reduced by a factor 2 as one can see in Fig.~\ref{bittner:fig_1}, where
the solution for a fixed $\tau$ is shown.
Obviously, for the low frequency signals towards the end of the divider chain a much sparser grid would be sufficient
for an accurate representation, in comparison to the high frequency input signal. 

\section{Division into subcircuits}

Although the representation of each signal over its own individual grid seems to give maximal flexibility, this approach
leads to several problems, which make the simulation inefficient. One problem is that the evaluation of the circuit
depends on the evaluation of device models, which is usually very costly, e.g.,\  for transistor models used nowadays.
Usually a device model has to be evaluated for every grid point. For a four terminal
 transistor we need therefore to evaluate
the device model for four different grids. In many cases this will be more costly than the evaluations for one optimized single
grid, which is against our intention to reduce the computational effort. This effort might be reduced if one has a strategy
to ``synchronize'' the grids, but that does not seem to be a trivial task.
On the other side, the signals show often similar
behavior, at least on parts of the circuit, such that the same grid might be (nearly) optimal for many signals. 
Thus, it might
be a better idea to collect signals of similar shape and use the same grid for all of these signals. Then we have to store only
a few grids, which makes it also easier to design an effective grid adaptation strategy.
Therefore, we consider groups of signals with similar
shape appearing in a part of the circuit. In the current implementation a priori knowledge
of the circuit design is required.

\begin{figure}
\sidecaption
\includegraphics[width=0.64\columnwidth]{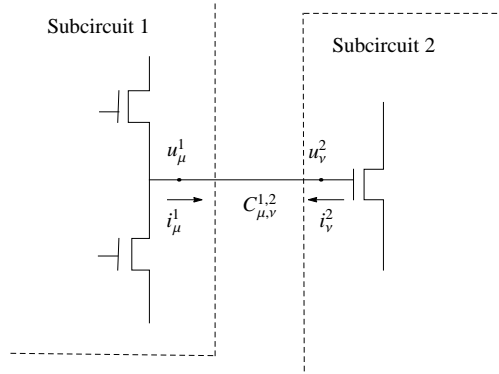}
\caption{Splitting of a circuit into subcircuits witch connections.}
\label{bittner:fig_subcircuit}       
\end{figure}

The circuit is divided into $N$ subcircuits which are connected at terminal nodes. To facilitate different expansions of signals
on the  subcircuits we replace each common node by a pair of nodes connected
by a perfect conductor, which is referred to as node tearing.
Namely, we introduce the ``connection'' $C_{\mu,\nu}^{k,\ell}$, if the $\mu$-th node of subcircuit $k$ is identified 
with the $\nu$-th node in subcircuit $\ell$, as one can see in Fig.~\ref{bittner:fig_subcircuit}. Thus, the circuit equations from
the modified nodal analysis (MNA) of the subcircuits have to be supplemented by additional conditions for the connections. 
The perfect conductor for the connection is modeled as voltage source of voltage zero. That is, we need the current through
the connection $C_{\mu,\nu}^{k,\ell}$, as additional unknowns $i_\mu^k$ and  $i_\nu^\ell$ for each of the two involved 
subcircuits. 

In addition to the resulting circuit equations
\begin{equation}\label{subcircuitmna}
\tfrac{d}{dt}q^k\big(x^k(t)\big) + g^k\big(x^k(t),t\big) =0,\qquad k=1,\ldots,N
\end{equation}
of the $N$ subcircuits, we need for each connection  $C_{\mu,\nu}^{k,\ell}$ that voltages and currents coincide,
that is we include the equations
\begin{equation}\label{connectionmna}
u_\mu^k(t)-u_\nu^\ell(t)=0 \mbox{\quad  and\quad} i_\mu^k(t)+i_\nu^\ell(t)=0.
\end{equation}
For the correct  understanding of the above formulation one needs to recall that $u_\mu^k(t)$ and $i_\mu^k(t)$ are
components of the vector $x^k(t)$ of unknowns, which contains \emph{all node voltages} (except ground) and currents 
through voltage sources, inductors, and \emph{connections}.

The splitting into subcircuits introduces the additional equations (\ref{connectionmna}), which will increase the problem
size. This may lead to a loss of performance if the splitting is done poorly. For a successful use of our method the splitting,
either done by hand or automatically, should follow some rules. First, a splitting should only be done if the signal shapes
in the resulting subcircuits differ enough to justify the use of different discretization grids, which are significantly coarser
than the grid for the (sub)circuit, which is splitted. Furthermore, the splitting should only generate few new \emph{connections}.
We expect that the second requirement is often fulfilled if the first requirement is satisfied  

\section{Spline Galerkin discretization and wavelet based adaptivity}

Our goal is to approximate the solution of the equations  (\ref{subcircuitmna})
and (\ref{connectionmna})  by spline functions as it was done in \cite{BiBra14b}. 
However, we want to use an adapted spline representation 
for each subcircuit, i.e.,
$$
x^i(t) = \sum_{k=1}^{n_i} c^i_k \varphi^i_k(t),\quad i=1,\ldots,N,
$$
where the families $\{\phi^i_k:~k=1,\ldots,n_i\}$ are periodic B-spline bases for spline spaces over grids of spline knots 
$T^i:=\{t^i_k\in(0,P]:~k=1,\ldots,n_i\}$, which may be mutually different.  We use a Petrov-Galerkin discretization to 
obtain a system of nonlinear equations, which determines the coefficients $c^i_k$. In particular, we integrate the equations 
(\ref{subcircuitmna}) and (\ref{connectionmna}) over subintervals, i.e.,
$$
\int_{\tau^i_{\ell-1}}^{\tau^i_\ell} \tfrac{d}{dt}q^i\big(x^i(t)\big) + g^i\big(x^i(t),t\big)\,dt =0, \qquad \ell = 1,\ldots,n_i,
$$
for each subcircuit
and
\begin{align}
\label{con1}
\int_{\tau^i_{\ell-1}}^{\tau^i_\ell} u^i_\mu(t)- u^j_\nu(t)\,dt&=0,\qquad  \ell = 1,\ldots,n_i\\
\label{con2}
\int_{\tau^j_{\ell-1}}^{\tau^j_\ell}  i^i_\mu(t)+ i^j_\nu(t)\,dt&=0,\qquad  \ell = 1,\ldots,n_j
\end{align}
for each connection $C^{i,j}_{\mu,\nu}$ between two subcircuits. 
The splitting points $\tau^i_\ell$ are chosen in close relation
to the spline grid, namely such the $t^i_\ell\in (\tau^i_{\ell-1},\tau^i_\ell)$. By using the grid $T^i$ in (\ref{con1}) but $T^j$ in (\ref{con2}), we assure that the number of unknowns and equations coincide.

The wavelet based coarsening and refinement methods described in \cite{BiBra14a,BiBra14b} are used to generate adaptive
grids for an efficient signal representation. An advantage of this approach is that grid and solution from the previous envelope
time step are used to generate an initial guess for Newton's method. Since the waveforms  change only slowly with $\tau$,
we have usually a very good initial guess on a nearly optimal grid and the solution is obtained with only few iteration steps.

\section{Numerical test}

The algorithm was implemented in 	C++ and tested on a PLL with frequency divider. The solutions for a fixed $\tau_k$
are shown in Figure.~1.

\begin{figure}
\centering
\includegraphics[width=\columnwidth]{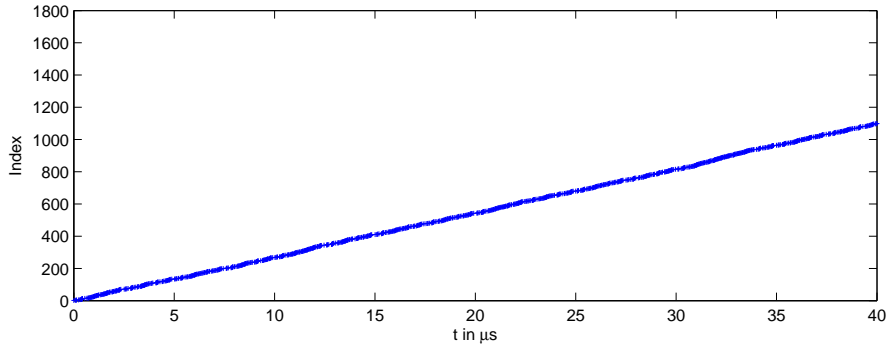}
\caption{Grid of the single grid method.}
\label{fig:bittner_2}       
\end{figure}

For comparison we show in Fig.~\ref{fig:bittner_2} the spline grid generated by the classical spline wavelet method 
(see \cite{BiBra14b}) using the same grid for all signals.
We have plotted the grid points $t_i$ against their index $i$, which allows to recognize the local density of the grid. 

\begin{figure}
\centering
\includegraphics[width=\columnwidth]{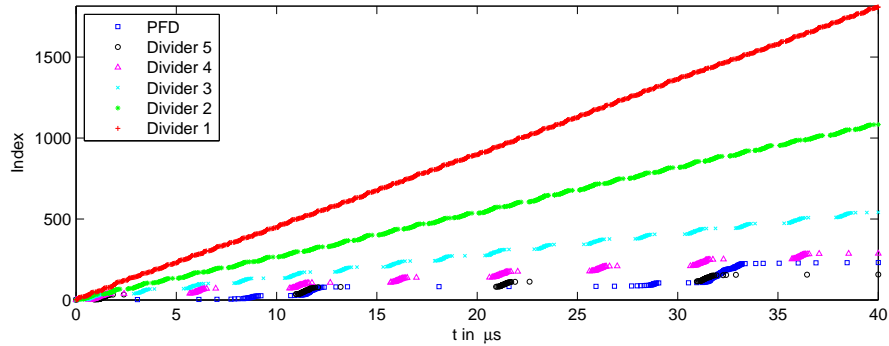}
\caption{Several signals in a frequency divider chain as part of a PLL.}
\label{fig:bittner_3}       
\end{figure}

The grids used in our new multiple grid method can be seen in Fig.~\ref{fig:bittner_3}. Obviously, one gets much better adapted, smaller grids for the lower frequency signals. This leads to a reduction of the total number of equations from roughly 130{,}000 to 85{,}000. The number of nonzeros in the Jacobian for Newton’s method is reduces from 5{,}000{,}000 to 2{,}500{,}000. Consequently the time for assembling resp. solving the linear system was reduced from 4s to 2s resp. 8s to 4s.

A further effect is that the larger the nonlinear system, the harder it is to solve by Newton’s method, which results in more Newton iterations as well as shorter envelope time steps. Thus, an envelope simulation with a frequency modulated signal over 0.3s worked well for the multiple grid method and was done in 37min. A similar simulation with the single grid method needed almost 5 hours.

\section{Conclusion}

An improvement of the spline wavelet based envelope method from \cite{BiBra14b} has been developed. It uses
different spline grids for different parts of the circuit. This leads to a more efficient representation of the solution, which results
in a significant reduction of computation time

\begin{acknowledgement}
This work has been partly supported by the ENIAC research project ARTEMOS under grant 829397 and
the FWF under grant P22549.
\end{acknowledgement}
%
%



\end{document}